\documentclass[12pt,a4paper]{amsart}

\usepackage{amsmath, amssymb, amsthm}
\usepackage{enumerate}

\makeatletter
\@namedef{subjclassname@2010}{%
  \textup{2010} Mathematics Subject Classification}
\makeatother

\newcommand{\bC}{\mathbb{C}}

\newcommand{\sP}{\mathsf{P}}
\newcommand{\bP}{\mathbb{P}}
\newcommand{\bN}{\mathbb{N}}
\newcommand{\rd}{\mathrm{d}}
\newcommand{\can}{\operatorname{can}}
\newcommand{\Id}{\mathrm{Id}}
\newcommand{\rep}{\operatorname{rep}}
\newcommand{\wan}{\operatorname{wan}}
\newcommand{\Lyap}{\operatorname{Lyap}}
\newcommand{\cS}{\mathcal{S}}
\newcommand{\supp}{\operatorname{supp}}
\newcommand{\sF}{\mathsf{F}}
\newcommand{\sJ}{\mathsf{J}}
\newcommand{\diam}{\operatorname{diam}}

\theoremstyle{plain}
\newtheorem{theorem}{Theorem}[section]
\newtheorem{lemma}[theorem]{Lemma}
\newtheorem{mainth}{Theorem}

\theoremstyle{definition}

\newtheorem*{acknowledgement}{Acknowledgement}
\theoremstyle{remark}

\renewcommand{\bold}[1]{\smallskip \noindent {\bf \boldmath #1 }\nopagebreak[4]}

\numberwithin{equation}{section}

\begin{document} 

\title[Approximation of Lyapunov exponents]{Approximation of Lyapunov exponents in non-archimedean and complex dynamics}
\author[Y\^usuke Okuyama]{Y\^usuke Okuyama}
\thanks{Partially supported by JSPS Grant-in-Aid for Young Scientists (B), 21740096 and 24740087.}
\address{
Division of Mathematics,
Kyoto Institute of Technology, Sakyo-ku,
Kyoto 606-8585 Japan}

\email{okuyama@kit.ac.jp}

\date{\today}

\begin{abstract}
We give two kinds of approximation of Lyapunov exponents 
of rational functions of degree more than one on the projective line
over more general fields than that of complex numbers.
\end{abstract}


\maketitle

\section{Introduction}\label{sec:intro}
Let $K$ be an algebraically closed field 
complete with respect to a non-trivial absolute value
$|\cdot|$. The field $K$ is said to be non-archimedean
if the strong triangle inequality $|z-w|\le\max\{|z|,|w|\}$ $(z,w\in K)$ 
holds (e.g.\ $p$-adic $\bC_p$). 
Otherwise, $K$ is archimedean, and then indeed $K\cong\bC$.
We always assume that $K$ has characteristic $0$.
The Berkovich projective line $\sP^1=\sP^1(K)$ produces a compactification of 
the (classical) projective line $\bP^1=\bP^1(K)$. 
For archimedean $K$, $\sP^1$ and $\bP^1$ are identical.
For the details of potential theory and dynamics on $\sP^1$, 
see \cite{BR10}, \cite{FR09}, \cite{Jonsson12}. 

Let $f$ be a rational function on $\bP^1$ of degree $d>1$. 
The action of $f$ on $\bP^1$ extends to 
a continuous, open, surjective and (fiber-)discrete map on $\sP^1$, 
preserving $\bP^1$ and $\sP^1\setminus\bP^1$.
The exceptional set $E(f)$ of (the extended) $f$ is
the set of all points $a\in\bP^1$ such that $\#\bigcup_{k\in\bN}f^{-k}(a)<\infty$.

For a rational function (say, possibly moving target) $a$ on $\bP^1$ and each $k\in\bN$,
there are exactly $d^k+\deg a$ roots of the equation $f^k=a$ 
in $\bP^1$ counting their multiplicity, unless $f^k\not\equiv a$. 
Let $\delta_w$ be the Dirac measure at $w\in\sP^1$.
Let us consider the averaged distribution
\begin{gather*}
 \nu_k^a:=\frac{1}{d^k+\deg a}\sum_{w\in\bP^1:f^k(w)=a(w)}\delta_w
\end{gather*}
of roots of $f^k=a$ in $\bP^1$,
where the sum takes into account the multiplicity of each root.
Let $\Delta$ be the normalized Laplacian on $\sP^1$.
The equilibrium (or canonical) measure of $f$ on $\sP^1$ is the Radon probability measure
\begin{gather*}
 \mu_f:=\Delta g_f+\Omega_{\can}
\end{gather*}
on $\sP^1$, where the continuous function $g_f$ is the dynamical Green function of $f$ on $\sP^1$
(cf.\ \cite[\S 2]{OkuFekete}),
and $\Omega_{\can}$ denotes the normalized Fubini-Study area element on $\bP^1$ for archimedean $K$,
and the Dirac measure at the canonical (or Gauss) point $\cS_{\can}$ on $\sP^1$ for
non-archimedean $K$. 

The equidistribution theorem for possibly moving targets,
which is in complex dynamics by Lyubich \cite[Theorem 3]{Lyubich83}
and in non-archimedean dynamics by Favre and Rivera-Letelier
\cite[Theor\`emes A et B]{FR09}, is

\begin{theorem}[{\cite[Theorem 3]{Lyubich83}, \cite[Th\'eor\`emes A et B]{FR09}}]\label{th:lyubich}
 Let $f$ be a rational function on $\bP^1=\bP^1(K)$ of degree $>1$.
 If a rational function $a$ on $\bP^1$ does not identically equal a
 value in $E(f)$, then $\nu_k^a\to\mu_f$ weakly on $\sP^1$ as $k\to\infty$.
\end{theorem}

Let $f^\#$ be the chordal derivative of $f$ with respect to
the normalized chordal distance $[\cdot,\cdot]$ on $\bP^1$. Then the function $f^\#$
extends continuously to $\sP^1$.
We define the Lyapunov exponent of $f$ with respect to $\mu_f$ by
\begin{gather*}
 L(f):=\int_{\sP^1}\log(f^\#)\rd\mu_f>-\infty
\end{gather*}
(cf.\ \cite[formula (1.2)]{OkuLog}). The (classical) critical set $C(f)$ of $f$ is 
\begin{gather*}
 C(f):=\{c\in\bP^1;f^\#(c)=0\}.
\end{gather*}
A point in the orbits of some periodic $c\in C(f)$ under $f$ is called a
superattracting periodic point of $f$.

We give two kinds of approximation of $L(f)$.

\bold{Repelling periodic points.}
A periodic point $p\in\bP^1$ of $f$ of period $k\in\bN$
is said to be repelling 
if $(f^k)^\#(p)=|(f^k)'(p)|>1$. For each $k\in\bN$, set
\begin{gather*}
 R_k(f):=\{p\in\bP^1;\text{ repelling periodic points of }f\text{ of period }k\},\\
 \nu_k^{\rep}:=\frac{1}{d^k+1}\sum_{w\in R_k(f)}\delta_w,\\
 R_k^*(f):=R_k(f)\setminus\bigcup_{j\in\bN;j<k,j|k}R_j(f),\quad
 \nu_k^*:=\frac{1}{d^k+1}\sum_{w\in R_k^*(f)}\delta_w.
\end{gather*}
From an argument based on Bezout's theorem, for every $k\in\bN$,
\begin{gather}
\#\left(\bigcup_{j\in\bN;j<k,j|k}R_j(f)\right)\le 2kd^{k/2}=o(d^k)\label{eq:Bezout}
\end{gather}
as $k\to\infty$ (cf.\ \cite[\S 4.2]{BDM08}).
If there are at most finitely many non-repelling periodic points of $f$ in $\bP^1$,
then Theorem \ref{th:lyubich} together with \eqref{eq:Bezout} implies 
the equidistribution
\begin{gather}
 \nu_k^{\rep}\to\mu_f\quad\text{and}\quad\nu_k^*\to\mu_f\label{eq:repconv} 
\end{gather}
weakly on $\sP^1$ as $k\to\infty$. In particular, for archimedean $K$,
from Fatou's finiteness on non-repelling periodic points of $f$,
\eqref{eq:repconv} always holds.

The following is a consequence of \cite[Theorem 1]{OkuLog}, which is a generalization of
\cite{BDM08} for archimedean fields (see also \cite{BertelootLyapunov}) and
of \cite{ST05} for number fields or function fields to general $K$.
The finiteness assumption is vacuous for archimedean $K$:

\begin{mainth}\label{th:superatt}
Let $f$ be a rational function on $\bP^1=\bP^1(K)$ of degree $d>1$. 
If there are at most finitely many non-repelling
periodic points of $f$ in $\bP^1$, then
\begin{gather}
\lim_{k\to\infty}\frac{1}{d^k+1}\sum_{w\in R_k(f)}\frac{1}{k}\log(f^k)^\#(w)=L(f),\label{eq:nonstrict}\\
\lim_{k\to\infty}\frac{1}{d^k+1}\sum_{w\in R_k^*(f)}\frac{1}{k}\log(f^k)^\#(w)=L(f),\label{eq:strict}
\end{gather}
and $L(f)\ge 0$.
\end{mainth}

Since this important case could be shown by 
a simpler argument than that in \cite[\S 4]{OkuLog}
and \eqref{eq:strict} is not mentioned there,
we give herewith a (simpler) proof of Theorem \ref{th:superatt}. 
The proof of Theorem \ref{th:superatt} for archimedean $K$ in \cite{BDM08}
was also based on the fact that $L(f)(\ge\log\sqrt{d})>0$ holds for archimedean $K$
(but this does not always hold for non-archimedean $K$).

\bold{Preimages of points.}
Put $B[z,r]:=\{w\in\bP^1;[w,z]<r\}$ for each $z\in\bP^1$ and each $r>0$.
Under the action $f$ on $\bP^1$, a point $z_0\in\bP^1$ is said to be {\itshape wandering}
if $\#\{f^k(z_0);k\in\bN\cup\{0\}\}=\infty$. 
Let $C(f)_{\wan}$ be the set of all $c\in C(f)$ wandering under $f$. The subset
\begin{gather*}
 E_{\wan}^{1/2}(f):=\bigcup_{c\in C(f)_{\wan}}
 \bigcap_{N\in\bN}\bigcup_{j\ge N}B[f^j(c),\exp(-d^{j/2})]
\end{gather*}
is of finite Hyllengren measure for the increasing sequence $(d^{j/2})\subset\bN$, so 
of (chordal) capacity $0$ (for a direct proof, see \cite[Lemma 2.1]{OkuFekete}). Set
\begin{gather*}
 E_{\Lyap}(f):=\left\{a\in\bP^1;\int_{\sP^1}\log(f^\#)\rd\nu_k^a\not\to L(f)\text{ as }k\to\infty\right\}.
\end{gather*}
This contains the finite set $\{f^j(c);c\in C(f)\setminus C(f)_{\wan},j\in\bN\}$
of all orbits of non-wandering (or preperiodic) critical points of $f$.
 
\begin{mainth}\label{th:nearlyeverywhere} 
 Let $f$ be a rational function on $\bP^1=\bP^1(K)$ of degree $d>1$. Then
\begin{gather*}
 E_{\Lyap}(f)\subset E_{\wan}^{1/2}(f)\cup\{f^j(c);c\in C(f)\setminus C(f)_{\wan},j\in\bN\}.
\end{gather*}
 In particular, for every $a\in\bP^1$ except for a set of capacity $0$,
 \begin{gather*}
  \lim_{k\to\infty}\frac{1}{d^k}\sum_{w\in f^{-k}(a)}\log(f^\#)(w)=L(f),
 \end{gather*}
 where the sum takes into account the multiplicity of each root $w$ 
 of the equation $f(\cdot)=a$ on $\bP^1$.
\end{mainth}

\section{Background}\label{eq:back}
Let $f$ be a rational function on $\bP^1=\bP^1(K)$ of degree $d>1$.
We fix a projective coordinate on $\bP^1$ so that $K\cong\bP^1\setminus\{\infty\}$.
Let $[\cdot,\cdot]$ be the chordal distance on $\bP^1$ normalized as
$[0,\infty]=1$ (for the precise definition, see, for example,
\cite[\S 2]{OkuLog}). 

For non-archimedean $K$,
a typical element $\cS\in\sP^1=\sP^1(K)$ other than $\infty$ is
regarded as a ($K$-closed) disk $\{z\in K;|z-a|\le r\}$ for some $a\in K$ and 
some $r=:\diam(\cS)\ge 0$.
The point $\cS_{\can}:=\{z\in K;|z|\le 1\}\in\sP^1$ is called
the canonical (or Gauss) point in $\sP^1$.
For disk $\cS$, put $|\cS|:=\sup_{z\in\cS}|z|$.
For disks $\cS,\cS'\in\sP^1$, let $\cS\wedge\cS'$ be the
smallest disk containing $\cS\cup\cS'$.

For non-archimedean $K$, $[\cdot,\cdot]$ canonically extends to
the generalized Hsia kernel $[\cS,\cS']_{\can}$ on $\sP^1$ with respect to
$\cS_{\can}$ satisfying 
for example that for disks $\cS,\cS'\in\sP^1$,
\begin{gather*}
 [\cS,\cS']_{\can}=\frac{\diam(\cS\wedge\cS')}{\max\{1,|\cS|\}\max\{1,|\cS'|\}}.
\end{gather*}
This extension is separately continuous on each variable $\cS,\cS'\in\sP^1$, and
vanishes if and only if $\cS=\cS'\in\bP^1$.
Let us also denote this extension by the same $[\cdot,\cdot]$,
for simplicity. 

For general $K$, the notion of (chordal) capacity of a Borel set in $\sP^1$
is introduced as usual using the chordal kernel $\log[\cdot,\cdot]$ on $\sP^1$.
We note that a countable set in $\bP^1$ is of capacity $0$.
The chordal potential of a Radon measure $\mu$ on $\sP^1$ is defined as
\begin{gather*}
 U_{\mu}^{\#}(\cdot):=\int_{\sP^1}\log[\cdot,\cS]\rd\mu(\cS)
\end{gather*}
on $\sP^1$. We recall the following lemmas.

\begin{lemma}[{\cite[Lemma 3.1]{OkuLog}}]\label{th:condition}
Let $f$ be a rational function on $\bP^1=\bP^1(K)$ of degree $d>1$. 
Suppose that a sequence $(\nu_k)$ of positive Radon measures on $\sP^1$
tends to $\mu_f$ weakly on $\sP^1$ as $k\to\infty$. Then
\begin{gather*}
 \lim_{k\to\infty}\int_{\sP^1}\log f^\#\rd\nu_k=L(f) 
\end{gather*}
holds if for each $c\in C(f)$,
$\lim_{k\to\infty}U_{\nu_k}^\#(c)=U_{\mu_f}^\#(c)$.
\end{lemma}

\begin{lemma}[{\cite[Lemma 3.3]{OkuLog}}]\label{th:potential}
Let $f$ be a rational function on $\bP^1=\bP^1(K)$ of degree $d>1$.
Let $a$ be a rational function on $\bP^1$ which does not identically
equal a value in $E(f)$, and 
let $(S_k)$ be a sequence of subsets $S_k\subset\sP^1$.
Then for every $z\in\bP^1\setminus(\limsup_{k\to\infty}\supp\nu_k^a)$, 
\begin{gather*}
 U_{\nu_k^a|(\sP^1\setminus S_k)}^\#(z)-U_{\mu_f}^\#(z)
=\frac{1}{d^k+\deg a}\log[f^k(z),a(z)]-U_{\nu_k^a|S_k}^\#(z)+o(1)
\end{gather*}
as $k\to\infty$.
\end{lemma}

The Berkovich Julia set $\sJ(f)$ is defined by the set of all
$\cS\in\sP^1$ satisfying
\begin{gather*}
 \bigcap_{U:\text{ an open neighborhood of }\cS\text{ in }\sP^1}\left(\bigcup_{k\in\bN}f^k(U)\right)=\sP^1\setminus E(f)
\end{gather*}
(\cite[Definition 2.8]{FR09}). 
The Berkovich Fatou set $\sF(f)$ is $\sP^1\setminus \sJ(f)$, which is open in $\sP^1$. 
We note that any superattracting (resp.\ repelling)
periodic points of $f$ is in $\sF(f)$ (resp.\ in $\sJ(f)$).

\section{A proof of Theorem \ref{th:superatt}}\label{sec:formula}
Set $a=\Id_{\bP^1}$ and 
$S_k=\{w\in\bP^1;f^k(w)=w\}\setminus R_k(f)$ for each $k\in\bN$.
Then $\nu_k^{\rep}=\nu_k^a|(\sP^1\setminus S_k)$. 
As seen in Introduction, 
under the assumptions in Theorem \ref{th:superatt},
the equidistribution \eqref{eq:repconv} holds.
In particular, 
\begin{gather}
 \lim_{k\to\infty}\nu_k^a|S_k=0\label{eq:degenerate}
\end{gather}
weakly on $\sP^1$.

Let us show that $\lim_{k\to\infty}U_{\nu_k^{\rep}}^\#(c)=U_{\mu_f}^\#(c)$ for each $c\in C(f)$.
Then Lemma \ref{th:condition} will conclude that
\begin{gather}
 \lim_{k\to\infty}\int_{\sP^1}\log(f^\#)\rd\nu_k^{\rep}=L(f),\label{eq:original}
\end{gather}
which with the chain rule shows \eqref{eq:nonstrict}. 

Let $c\in C(f)\cap\sJ(f)$. Then $c$ is not periodic under $f$, or equivalently,
$c\in\bP^1\setminus(\bigcup_{k\in\bN}\supp\nu_k^a)$. 
In particular, $c\not\in\bigcup_{k\in\bN}S_k$, and
under the assumption in Theorem \ref{th:superatt}, which is equivalent to
$\#\bigcup_{k\in\bN}S_k<\infty$,
we have $\inf\{[c,w];w\in \bigcup_{k\in\bN}S_k\}>0$. 
This with \eqref{eq:degenerate} implies that 
\begin{gather*}
 \lim_{k\to\infty}U_{\nu_k^a|S_k}^\#(c)=0.
\end{gather*}
Moreover, Przytycki \cite[Lemma 1]{Przytycki93} asserts that 
{\itshape for any $c\in C(f)\cap\sJ(f)$ and every $k\in\bN$,
\begin{gather*}
 [f^k(c),c]\ge \frac{1}{10}(\max\{1,\sup_{\sP^1}f^{\#}\})^{-k+1}
\end{gather*}}
(the original proof of \cite[Lemma 1]{Przytycki93} for archimedean $K$
works for non-archimedean $K$), so for every $c\in C(f)\cap\sJ(f)$,
\begin{gather*}
 \lim_{k\to\infty}\frac{1}{d^k+1}\log[f^k(c),c]=0. 
\end{gather*}
Hence Lemma \ref{th:potential} implies that $\lim_{k\to\infty}U_{\nu_k^{\rep}}^\#(c)=U_{\mu_f}^\#(c)$.

Next, let $c\in\sF(f)$ (and we will not use Lemma \ref{th:potential}). 
Then $\log[c,\cdot]$ is continuous, i.e, does not take $-\infty$,
on $\sJ(f)$. Since
$\bigcup_{k\in\bN}R_k(f)\subset\sJ(f)$, the equidistribution
\eqref{eq:repconv} implies that $\lim_{k\to\infty}U_{\nu_k^{\rep}}^\#(c)=U_{\mu_f}^\#(c)$. 

Now the proof of \eqref{eq:original} is complete.

Finally, from $1\le\inf_{\bigcup_k R_k(f)}f^\#\le\sup_{\bP^1}f^\#<\infty$ and \eqref{eq:Bezout},
we have
\begin{gather*}
 \lim_{k\to\infty}\int_{\sP^1}\log
 f^\#\rd\nu_k^*=\lim_{k\to\infty}\int_{\sP^1}\log f^\#\rd\nu_k^{\rep}
\end{gather*}
(if one of the limits exists), 
which together with \eqref{eq:original} and the chain rule shows \eqref{eq:strict}. 
Now the last assertion $L(f)\ge 0$ is obvious.\qed

\section{A proof of Theorem \ref{th:nearlyeverywhere}}

We set $S_k=\emptyset$ for each $k\in\bN$.
Let $a\in\bP^1\setminus (E_{\wan}^{1/2}(f)\cup\{f^j(c);c\in C(f)\setminus C(f)_{\wan},j\in\bN\})$.
Then for every $c\in C(f)$, we have $c\in\bP^1\setminus(\limsup_{k\to\infty}\supp\nu_k^a)$ and
\begin{gather*}
 \liminf_{k\to\infty}\frac{1}{d^k}\log[f^k(c),a]\ge\liminf_{k\to\infty}(-d^{-k/2})=0,
\end{gather*}
which with $[\cdot,\cdot]\le 1$ implies that $\lim_{k\to\infty}(\log[f^k(c),a])/d^k=0$. 
From this and the assumption that $S_k=\emptyset$ for each $k\in\bN$,
Lemma \ref{th:potential} implies that $\lim_{k\to\infty}U_{\nu_k^a}^{\#}(c)=U_{\mu_f}^\#(c)$.
From this,
Lemma \ref{th:condition} completes the proof of Theorem \ref{th:nearlyeverywhere}. \qed

\begin{acknowledgement}
 The author thanks the referee for useful comments.
\end{acknowledgement}
 
\def\cprime{$'$}

\end{document}